\documentclass[a4paper,10pt]{article}
\usepackage[utf8x]{inputenc}
\usepackage{amssymb} 

\title{On a complexification of the moduli space of Bohr - Sommerfeld lagrangian cycles}
\author{Nikolai A. Tyurin\footnote{The author is partially supported by Laboratory of Mirror Symmetry NRU HSE, RF Government grant, ag. N 14.641.31.0001} \\  Lab of Theoretical  Physics, JINR (Dubna) and\\
 Lab of Mirror Symmetry, NRU HSE (Moscow) }

\begin{document}

\maketitle

\begin{abstract} In the previous papers we present a construction of the set ${\cal U}_{SBS}$ in the direct product ${\cal B}_S \times \mathbb{P} \Gamma (M, L)$
of the moduli space of Bohr - Sommerfeld lagrangian submanifolds of fixed topological type and the projectivized space of smooth sections of the prequantization bundle
$L \to M$ over a given compact simply connected symplectic manifold $M$. Canonical projections $p: {\cal U}_{SBS} \to \mathbb{P} \Gamma (M, L)$ and $q: {\cal U}_{SBS}
\to {\cal B}_S$ are studied in the present text: first, we show that the differential ${\rm d} p$ at a given point is an isomorphism, which implies that a natural complex structure
can be defined on ${\cal U}_{SBS}$; second, the projection $q: {\cal U}_{SBS} \to {\cal B}_S$ splits as the combination ${\cal U}_{SBS} \to T {\cal B}_S \to {\cal B}_S$
such that the fibers of the first map are complex subsets in ${\cal U}_{SBS}$. This implies that an appropriate section of the first map should define
a  complex structure on $T {\cal B}_S$; therefore it can be seen as a complexification of the moduli space ${\cal B}_S$. The construction can be exploited
in the Lagrangian approach to Geometric Quantization. 

\end{abstract}

\section{Introduction}

Let $(M, \omega)$ be a compact simply connected symplectic manifold of dimension $2n$ with integer symplectic form so $[\omega] \in H^2(M, \mathbb{Z}) \subset H^2 (M, \mathbb{R})$.
Consider tha prequantization data $(L, a)$ where $L \to M$ is a complex line bundle with fixed hermitian structure $h$ and $a \in {\cal A}_h(L)$ is a hermitian connection
such that the curvature $F_a = 2 \pi i \omega$. Fixing a topological type ${\rm top} S$ of a smooth orientable $n$ - dimensional manifold and homology class
$[S] \in H_n(M, \mathbb{Z})$ one can consider the moduli space ${\cal B}_S$ of Bohr - Sommerfeld lagrangian cycles of fixed topological type, constructed in [1].
This moduli space is an infinite dimensional Fr\'{e}chet smooth real manifold, locally modelled by the spaces $C^{\infty}(S, \mathbb{R})$ modulo constants. Points of
${\cal B}_S$ can be described as lagrangian submanifolds of the fixed topological type which satisfy the Bohr - Sommerfeld condition: the restriction of the prequantization data
$(L, a)|_S$ admits covariantly constant sections. The details can be found in [1].

Take the space   $\Gamma(M, L)$ of all smooth sections of the prequantization bundle and consider its projectivization $\mathbb{P} \Gamma(M, L)$. Every point of the last projective space
corresponds to a complex 1 - form defined on an open subset $M \backslash D_{\alpha}$. Namely for any point $[\alpha] \in \mathbb{P} \Gamma (M, L)$ take a smooth section $\alpha$ which is correctly
defined up to $\mathbb{C}^*$ and consider the following 1 -form
$$
\rho(\alpha) = \frac{\nabla_a \alpha}{\alpha} = \frac{<\nabla_a \alpha, \alpha>_h}{<\alpha, \alpha>_h} \in \Omega^1_{M \backslash D_{\alpha}} \otimes \mathbb{C},
$$
where $D_{\alpha}$ is the zeroset of $\alpha$. It is clear that $\rho(\alpha)$ is the same for any rescaled $c \alpha$ thus it corresponds to $[\alpha] \in \mathbb{P} \Gamma(M, L)$.

The main properties of $\rho(\alpha)$ are the following: its real part ${\rm Re} \rho(\alpha)$ is exact, and the differential of the imaginary part satisfies ${\rm d} {\rm Im} \rho(\alpha) = 2 \pi \omega$
(see, f.e. [2]). 

Then we say that a lagrangian submanifold $S \subset M$ is special Bohr - Sommerfeld with respect to $\alpha$ (or $\alpha$ - SBS for short) iff the restriction ${\rm Im} \rho(\alpha)|_S$ identically vanishes.
In particular this means that $S$ does not intersect the zeroset $D_{\alpha}$. Thus the SBS - condition derives a subset ${\cal U}_{SBS} \subset {\cal B}_S \times \mathbb{P} \Gamma(M, L)$ formed by pairs $(S, [\alpha])$ such that $S$ is $\alpha$ - SBS.

{\bf Remark.} In paper [2] we deduced the last defintion from the first one, which we omitt here for breivity. 

By the very definition the subset ${\cal U}_{SBS}$ admits two canonical projections $p: {\cal U}_{SBS} \to \mathbb{P} \Gamma (M, L)$ and $q: {\cal U}_{SBS} \to {\cal B}_S$.
The first projection has been studied in [2], where one established that the fibers of the projections are discrete; the image of the projection is an open set in the projective space;
no ramification takes place; the differential $d p$ has trivial kernel at generic point. From this one has deduced that ${\cal U}_{SBS}$ is  weakly Kahler: the standard Kahler structure
on the projective space $\mathbb{P} \Gamma (M, L)$ can be lifted there. In the finite dimensional case the triviality of the kernel would imply that ${\rm d} p$ is an isomorphism, but
the infinite dimensional situation is much more complicated. There are examples when the kernel of the map is trivial but the map itself is not an isomorphism. 

Below we use the correspondence $[\alpha] \longleftrightarrow \rho(\alpha)$ to prove that actially the differential ${\rm d} p$ is an isomorphism: any tangent vector to $\mathbb{P} \Gamma(M, L)$
can be canonically lifted to ${\cal U}_{SBS}$; in particular a complex structure on ${\cal U}_{SBS}$ can be presented explicitly in terms of the lifting.

On the other hand the space ${\cal U}_{SBS}$ can not be considered as a compexification of the moduli space ${\cal B}_S$: it is too big.   Indeed, the  ``dimension'' of ${\cal B}_S$ is 
$C^{\infty}(S, \mathbb{R})$ modulo constants while the  ``dimension'' of ${\cal U}_{SBS}$ is the same as of $\mathbb{P} \Gamma(M, L)$ and the last one is $C^{\infty}(M, \mathbb{C})$ modulo
constants. At the same time it exists a natural map $\tau: {\cal U}_{SBS} \to T {\cal B}_S$ which factorizes the first projection $q$, namely $q = \tau \circ \pi$, where $\pi: T {\cal B}_S
\to {\cal B}_S$ is the canonical projection. Tha map $\tau$ is defined very naturally: since by the very definition for any point $(S, [\alpha]) \in {\cal U}_{SBS}$
one has that $\rho(\alpha)|_S$ is an exact real 1 -form and at the same time any real exact 1 -form is a tangent vector to ${\cal B}_S$ at the point $S$ we get simple formula 
$$
\tau(S, [\alpha]) = (\rho(\alpha))|_S \in T_S {\cal B}_S,
$$
and evidently applying $\pi$ we get the result of $q(S, [\alpha]) = S \in {\cal B}_S$.

Below we prove that a generic fiber $\tau^{-1}(S, {\rm d}f) \subset {\cal U}_{SBS}$ is a complex subset of ${\cal U}_{SBS}$ with respect to the complex structure lifted from the projective space $\mathbb{P}\Gamma(M, L)$.
This should lead to a construction of a complex structure on $T {\cal B}_S$: suppose one finds a natural section of the fibration $\tau: {\cal U}_{SBS} \to T {\cal B}_S$ which is Kahler orthogonal to the fibers
or is just a complex subset of ${\cal U}_{SBS}$. Right now we can not present a good candidate, hoping to find an approriate one in a future.

{\bf Acknowledgments.} This work was done as an extended comment to a short lecture course given at the VII school on Geometry and Physics in Bialowieza. The author cordially thanks  the organizers and the participants
of the shool for hospitality, valuable discussions and remarks. 

\section{``Affinitzation'' of the projective space $\mathbb{P} \Gamma (M, L)$}

To construct the complex structure on ${\cal U}_{SBS}$ via the projection $p$ first of all we study the correspondence $[\alpha] \leftrightarrow \rho(\alpha)$ where $\alpha \in \Gamma(M, L)$ is a smooth section
of the prequantization bundle and $\rho(\alpha)$ is the corresponding 1 - form given by the formula
$$
\rho(\alpha) = \frac{\nabla_a \alpha}{\alpha} \in \Omega^1_{M \backslash D_{\alpha}} \otimes \mathbb{C}.
$$
The main properties of the last 1 - form read as follows
$$
{\rm Re} \rho(\alpha) = {\rm d} (ln \vert \alpha \vert_h) - \textbf{exact}, \quad {\rm d} {\rm Im} \rho(\alpha) = 2 \pi 
$$
on the complement $M \backslash D_{\alpha}$, where $D_{\alpha}$ is the zeroset. Below we suppose that $\alpha$ is generic, in particular it means that $D_{\alpha}$ is not too big and pathologic. But in any case 
we have

{\bf Proposition 1.} {\it For any pair $\alpha_1, \alpha_2$ of smooth sections with zeros $D_{\alpha_1}, D_{\alpha_2}$ if $\rho(\alpha_1) \equiv \rho(\alpha_2)$
on $M \backslash \{D_{\alpha_1} \cup D_{\alpha_2}\}$ then there $\alpha_2 \equiv c \alpha_1$ for a complex constant $c \in \mathbb{C}$.}

Indeed, since $\alpha_i$ are smooth sections of a complex line bundle it exists a non vanishing complex function $\psi$ such that $\alpha_2 = \psi \alpha_1$
on $M \backslash \{D_{\alpha_1} \cup D_{\alpha_2}\}$. Calculating the corresponding $\rho$ - forms one gets
$$
\rho(\alpha_2) = \frac{\nabla_a \psi \alpha_1}{\psi \alpha_1}  = \frac{{\rm d} \psi}{\psi} + \frac{\nabla_a \alpha_1}{\alpha_1} \equiv \rho(\alpha_1)
\eqno (1)
$$
 by the assumption. Thus we deduce the identity $\frac{{\rm d} \psi}{\psi} \equiv 0$ on $M \backslash \{D_{\alpha_1} \cup D_{\alpha_2}\}$ which implies $\psi \equiv const$.

 Moreover, the $\rho$ - form defines the corresponding class almost uniquelly: it remains to discuss the role of the zeroset $D_{\alpha}$, namely if we have a pair $({\rm d} f + \imath \lambda, D_{\alpha})$
where ${\rm d} f + \imath \lambda$ is finite outside of $D_{\alpha} \subset M$, ${\rm d} \lambda = 2 \pi \omega$ and $D_{\alpha}$ is the zeroset of a smooth section $\alpha$ then it is possible to say
when ${\rm d} f + \imath \lambda$ equals to $\rho(\alpha')$ for certain $\alpha'$. The difference ${\rm d} f + \imath \lambda - \rho(\alpha)$ is always of the form ``exact form + i closed form'',
and for the existence of such an $\alpha'$ the last closed form in the difference must present an integer class from $H^1(M \backslash D_{\alpha}, \mathbb{R})$.

From this we see that while the space of classes $[\alpha]$ is projective the space of $\rho$ - forms looks much more linear: if we forget about zerosets then this space is affine
(the differences lie in the space of closed forms which present integer cohomology classes). At least this ``affinization'' works at the local level: for a class of smooth sections
$[\alpha]$ there is a small neighborhood ${\cal O}([\alpha]) \subset \mathbb{P} \Gamma (M, L)$ such that all classes from ${\cal O}([\alpha])$ are uniquelly determined
by pairs ${\rm d} f + \imath {\rm d} g \in \Omega^1_{M \backslash {\cal O}( D_{\alpha})} \otimes \mathbb{C}$. At the infinitesimal level: take a family $\alpha_t, t \in [0;1]$ and conisder 
the difference form $\Delta_t = \frac{1}{t}(\rho(\alpha_t) - \rho(\alpha_0))$. For sufficiently small $t$ this form equals ${\rm d} f_t + \imath {\rm d} g_t$, an exact form,  which tends to the form
${\rm d} f_0 + \imath {\rm d} g_0$, and the last one is defined on $M \backslash D_{\alpha_0}$ and understood as a tangent vector to the $\rho$ - space.

On the other hand for any ${\rm d} f_0 + \imath {\rm d} g_0$ on $M \backslash D_{\alpha_0}$ one can take 
$\alpha_t = e^{t(f_0 + \imath g_0)} \alpha_0$ realizing this tangent vector by the simplest  tajectory in $\mathbb{P} \Gamma (M, L)$ which is just a segment in the $\rho$ - space. 
Indeed, according to formula (1) above the difference  $\rho(\alpha_t) - \rho(\alpha_0)$ equals to $\frac{{\rm d} \psi}{\psi}  = t ({\rm d} f_0 + \imath {\rm d} g_0)$, therefore $\Delta_t \equiv {\rm d} f_0 
+ \imath {\rm d} g_0$ for this trajectory.

Thus we can consider the $\rho$ - representation for the projective space $\mathbb{P}\Gamma(M, L)$ since for our purposes it is much simpler to study pairs $(\rho(\alpha), D_{\alpha})$
instead of equivalence classes of smooth sections. It remains to present the corresponding expression for the complex structure on $\mathbb{P} \Gamma(M, L)$
in the $\rho$ - coordinates.

Recall, that our space $\Gamma(M, L)$ is naturally endowed with the hermitian scalar product $\int_M <\alpha_1, \alpha_2>_h d \mu_L$, where $d \mu_L$ is the Liouville volume form.
Therefore the tangent space $T_{[\alpha_0]} \mathbb{P} \Gamma(M, L)$ is given by sections $\{\delta \alpha \}$ such that $\int_M <\delta \alpha, \alpha_0>_h d \mu_L = 0$.
Then the standard complex structure acts by the multiplication $\delta \alpha \mapsto \imath \delta \alpha$, and we must find the correspondence between $\delta \alpha$
and $\Delta_t$ in the $\rho$ - representation. For this take the family of sections $\alpha_t = \alpha_0 + t \delta \alpha$ where $\delta \alpha$ is perpendicular to
$\alpha_0$ and consider $\Delta_t = \frac{1}{t}(\rho(\alpha_0 + t \delta \alpha) - \rho(\alpha_0))$. Simple calculation leads to
$\Delta_t = \frac{\nabla_a \delta \alpha}{\alpha_0} + o(t)$ therefore $\frac{\nabla_a \delta \alpha}{\alpha_0} = {\rm d} f_0 + \imath {\rm d} g_0$ for certain
$f_0, g_0$. But if we multiply $ \delta \alpha$ by $\imath$ then it just corresponds to the multiplication 
$$
\frac{\nabla_a \imath \delta \alpha}{\alpha_0} = \imath \frac{\nabla_a \delta \alpha}{\alpha_0} = {\rm d} g_0 - \imath {\rm d} f_0,
$$
it follows that the complex structure in the $\rho$ - representation acts just as the natural multiplication by $\imath$, which is possible since the symmetry: both
the real and imaginary parts are exact 1- forms.  

\section{Hamiltonian deformations}

Now we start to prove the main technical fact: the differential 
$${\rm d} p: T_{(S, [\alpha])} {\cal U}_{SBS} \to \mathbb{P} \Gamma(M, L)
$$
 is an isomorphism at generic point
$(S, [\alpha])$. In our arguments we will use the $\rho$ - representation of the projective space, therefore a point of ${\cal U}_{SBS}$ will be represented by 
pair $(S, \rho(\alpha))$. Recall that the elements of the pair satisfy $\rm{Im} \rho(\alpha)|_S \equiv 0$. 

At the first step we establish an important dynamical property of the subset ${\cal U}_{SBS}$: take any smooth function $F$ on $M$ and generate the Hamiltonian flow
$\phi^t_{X_F}$ where $X_F$ is the Hamiltonian vector field of $F$. Then

{\bf Proposition 2.} {\it If $(S_0, \rho_0) \in {\cal U}_{SBS}$ then $(\phi^t_{X_F}(S_0), (\phi^t_{X_F})^*(\rho_0)) \in {\cal U}_{SBS}$.}

Indeed, the condition $\rm{Im} \rho|_S \equiv 0$ is stable with respect to symplectomorphisms. In particular the infinitesimal part of this deformation
gives a vector field $\Theta(F) \in \rm{Vect} ({\cal U}_{SBS})$. The components of this vector field can be separated with respect to the elements of the pair:
for the Bohr - Sommerfeld lagrangian submanifold $S_0$ the tangent component equals to ${\rm d} F|_{S_0} \in T_{S_0} {\cal B}_S$; for the $\rho$ - component
the tangent component is given by the Lie derivative ${\cal L}_{X_F} \rho_0 = {\rm d} (\rm{Re} \rho_0 (X_F)) + \imath ({\rm d} (\rm{Im} \rho_0(X_F)) +  2 \pi {\rm d} F)$,
which is the sum of exact forms. Complete description must include a normal vector field on $D_{\alpha_0}$ since the flow $\phi^t_{X_F}$ deforms
the zeroset $D_{\alpha_0}$ as well but it does not affect the story near $S_0$ and we will skip it. 

The next step: study the situation near our Bohr - Sommerfeld submanifold $S_0$, namely in a sufficiently small Darboux - Weinstein neighborhood ${\cal O}_{DW}(S_0) 
\subset M$ which is by the very definition symplectomorphic to a small neighborhood of the zero section in the cotangent bundle $T^* S_0$. Using this symplectomorphism
we transport the canonical action form $\rho_{can}$ from the cotangent space to ${\cal O}_{DW}(S_0)$ and denote it by the same symbol $\rho_{can}$.
Then as it was shown in [2], the difference $\rho_{can} - \frac{1}{2 \pi} \rm{Im} \rho_0$ is an exact form ${\rm d} \Psi$ totaly vanishing at $S_0$.
This means that we can deform the Darboux - Weinstein neighborhood ${\cal O}_{DW}(S_0)$ such that the new canonical action form shall be equal
to $\rm{Im} \rho_0$. Indeed, for this we can find a Hamiltionian transformation generated by certain function $F$ which moves $\rho_{can}$ to
$\rho_{can} + {\rm d} \Psi$. This means that we need to solve the equation ${\cal L}_{X_F} \rho_{can} = {\rm d} \Psi$. In the local Darboux coordinates this equation 
looks as follows:
$$
\sum_{i=1}^n p_i \frac{\partial F}{\partial p_i} + F = \Psi;
\eqno (2)
$$
decomposing the left and right hand sides in the powers of $p_i$ and taking into the account that $\Psi = o(p_i)$ we get the desired $F$.
This $F$ vanishes  on $S_0$ therefore the corresponding flow does not move $S_0$. Making if necessury the neighborhood ${\cal O}_{DW}(S_0)$ smaller we get
at the end a Darboux - Weinstein neighborhood of $S_0$ such that $\rho_{can} \equiv \rm{Im} \rho_0|_{{\cal O}_{DW}(S_0)}$.

The next step: suppose we have a tangent vector $\delta \rho = {\rm d} f_0 + \imath {\rm d} g_0$ at point $\rho_0$ in the $\rho$ - space. Our goal is to lift it
to the tangent space $T_{(S_0, \rho_0)} {\cal U}_{SBS}$. Since in the defining property of ${\cal U}_{SBS}$ we do exploit the imaginary part
of $\rho$ only it is sufficient to find a hamiltonian deformation $\phi^t_{X_F}$ which realizes $\rho_0 \mapsto \rho_0 + \rm{something} + \imath {\rm d} g_0$,
where ${\rm d} g_0$ is a fixed tangent component for the imaginary part. Then for this Hamiltonian deformation we take the corresponding deformation of $S_0$
which we denote as $\delta S_0(\delta \rm{Im} \rho_0)$. Then the total deformation of the pair $(S_0, \rho_0)$ is given by $\Theta (F)$;
and it remains to mention that the variation of the real part of $\rho_0$ does not affect the defining condition $\rm{Im} \rho|_S \equiv 0$ therefore
$\delta \rho_0 = {\rm d} f_0 + \imath {\rm d} g_0$ is lifted to $(\delta S_0({\rm d} g_0), {\rm d} f_0 + \imath {\rm d} g_0)$. 

To realize this scheme we need to show that any deformation of the imaginary part of $\rho_0$ can be presented as the result of a Hamiltonian deformation
generated by certain function $F$. But   since we have rearranged the Darboux - Weinstein neighborhood such that $\rm{Im} \rho_0 = \rho_{can}$
the equation for this $F$ shall be essentially the same as above, namely the equation (2) with the fixed right hand side $g_0$ gives us the desired
$F$. Again using the decomposition in the powers of $p_i$ we establish the existance of such an $F$; at the same time the infinitesimal deformation of a given
Bohr - Sommerfeld submanifold corresponds to an exact 1 - form on it, and if the Hamiltonian deformation is generated by function $F$ then the corresponding 
tangent vector in $T_{S_0}{\cal B}_S$ equals to ${\rm d} F|_{S_0}$. From the equation (2) it is clear that the desired $F$ must satisfies $F|_{S_0} \equiv
g_0|_{S_0}$ which leads to the final form of the lifting:

{\bf Proposition 3.} {\it At  generic point $(S_0, \rho_0) \in {\cal U}_{SBS}$ in the $\rho$ - representation the inverse linear map $({\rm d} p)^{-1}$
is given by the formula $({\rm d} p)^{-1}({\rm d} f_0 + \imath {\rm d} g_0) = ({\rm d} (g_0|_{S_0}); {\rm d} f_0 + \imath {\rm d} g_0)$.}

In particular we have a lifted complex structure $\tilde I: T {\cal U}_{SBS} \to T {\cal U}_{SBS}$, acting as follows
$$
\tilde I ({\rm d} (g_0|_{S_0}), {\rm d} f_0 + \imath {\rm d} g_0) = ({\rm d} (f_0|_{S_0}), - {\rm d} g_0 + \imath {\rm d} f_0).
\eqno (3)
$$

From the arguments above we can see that the difficulties at non generic points of ${\cal U}_{SBS}$ can appear due to the behavier of the zerosets $D_{\alpha}$;
at the same time by the very definition for any pair $(S, [\alpha]) \in {\cal U}_{SBS}$ we have $S \cap D_{\alpha} = \emptyset$ therefore it is possible
to consider the Darboux - Weinstein neighborhoods with the same property ${\cal O}_{DW} (S) \cap D_{\alpha} = \emptyset$ and it seems that the proof must work
for any point. Anyway the last formula (3) is so simple that the operator $\tilde I$ can be extended from the generic points of ${\cal U}_{SBS}$ to the whole. 

\section{Fibers of the map $\tau: {\cal U}_{SBS} \to T {\cal B}_S$}

The geometry of the subset ${\cal U}_{SBS}$ is quite reach: there are many natural conditions which cut subsets, subspaces etc. 

For example choose any smooth function $f \in C^{\infty}(M, \mathbb{R})$ and define the following natural subset ${\cal U}^f \subset {\cal U}_{SBS}$ by the condition
$$
{\cal U}^f = \{  (S, \rho) \in {\cal U}_{SBS} \quad \vert \quad  \rho|_S \equiv {\rm d} f|_S \}.
$$

 The intersection of two such subsets for different functions $f_1, f_2$ has some special projection to ${\cal B}_S$: since by the very definition
$(S, \rho) \in {\cal U}^{f_1} \cap {\cal U}^{f_2}$ must satisfy ${\rm d} f_1|_S = {\rm d} f_2 |_2$ then the Bohr - Sommerfeld submanifold $S$ must be stationary
with respect to the Hamiltonian action of the difference $f_1 - f_2$. This means that $S$ lies in a level set of the function
$f_1 - f_2$.  Moreover if we take three functions $f_1, f_2, f_3$ then the triple intersection
${\cal U}^{f_1} \cap {\cal U}^{f_2} \cap {\cal U}^{f_3}$ satisfies even more sharp condition. Indeed, if $(S, \rho)$ lies in the triple intersection then
$S$ must lie in the intersection of two level sets of $f_3 - f_1$ and $f_2 -f_1$. Therefore despite of the size of each ${\cal U}^{f_i}$ the intersection
${\cal U}(f_1, ..., f_n) = {\cal U}^{f_0} \cap ... {\cal U}^{f_{n}}$ is empty for generic set of functions $f_1, ..., f_n$ where we choose $f_0 \equiv 0$. Exceptionly the cases when
$f_1, ..., f_n$ are not generic can be distinguished in this circumstance: if, say, $f_1, ..., f_n$ commute then the intersection is non empty, having as a``support''
Bohr - Sommerfeld fibers of the corresponding action map $F_{act} = (f_1, ..., f_n): M \to \mathbb{R}^n$. The number of such fibers is always finite, 
and they play important role in Geometric Quantization procedure in the presence on real polarization, see [3]. The Hilbert space for the corresponding
quantization is spanned by the finite number of Bohr - Sommerfeld lagrangian fibers ${\cal H} = \mathbb{C} <S_1> \oplus ... \oplus \mathbb{C} <S_m>$, see [3],
and in a sence our subset ${\cal U}(f_1, ..., f_n)$ has a finite limit  given by the projectivization $\mathbb{P} {\cal H}$. 

In [3] one formulates the main problem which appears in this approach, namely the problem of transition which should relate two quantum systems given
by different sets of commuting functions $(f_1, ..., f_n), (g_1, ..., g_n)$. Bohr - Sommerfeld lagrangian fibers for the corresponding action maps 
$F_{act}^f = (f_1, ..., f_n)$ and $F_{act}^g = (g_1, ..., g_n)$ are different, and it is possible to calculate the transition amplitudes for the corresponding Hilbert space
if two lagrangian fibrations are realted by a Hamiltonian deformation. For generic case it is not quite clear how to define the amplitudes. 

May be these amplitudes can be calculated via certain realtions between two subsets ${\cal U}(f_1, ..., f_n)$ and ${\cal U}(g_1, ..., g_n)$ or their projections
$p({\cal U}(f_1, ..., f_n))$ and $p({\cal U}(g_1, ..., g_n))$ in the projective space $ \mathbb{P} \Gamma (M, L)$.

There is another natural construction which gives a map $\tau: {\cal U}_{SBS} \to T {\cal B}_S$. The definition is extremely simple:  for any pair $(S_0, \rho_0= \rho(\alpha_0)) \in {\cal U}_{SBS}$
the restriction  $\rho_0|_{S_0}$ is an exact 1 -form, and any exact 1 -form represents a tangent vector from $T_{S_0} {\cal B}_S$. Thus
$$
\tau(S_0, \rho_0) = (S_0, \rho_0|_{S_0}) \in T {\cal B}_S.
$$

For this map we have the following

{\bf Proposition 4.} {\it Any fiber of the map $\tau: {\cal U}_{SBS} \to T {\cal B}_S$ is a complex subset in ${\cal U}_{SBS}$ with respect to the complex structure $\tilde I$, defined in (3).}

Choose a point $S_0 \in {\cal B}_S$ and a tangent vector $v_0 = ({\rm d} h_0|_{S_0}) \in T_{S_0} {\cal B}_S$ given by an exact 1 - form on $S_0$. Then the fiber $\tau^{-1}(S_0, v_0)$
can be described as follows: every $(S_0, \rho) \in \tau^{-1}(S_0, v_0)$ contains as the second element a form $\rho$ with fixed restriction to $S_0$ namely 
$\rm{Re} (\rho|_{S_0}) = {\rm d} h_0|_{S_0}, \rm{Im} (\rho|_{S_0}) \equiv 0$. Therefore the infinitesimal variations of such forms satisfies 
$$\delta \rho = {\rm d} f_0 + \imath {\rm d} g_0 \quad | \quad {\rm d} f_0|_{S_0} = {\rm d} g_0 |_{S_0} \equiv 0.
$$

This implies that the tangent space to the fiber is complex with respect to the complex structure $\tilde I$, which complites the proof of the Proposition.

It follows that appropriate sections of the fibration $\tau: {\cal U}_{SBS} \to T {\cal B}_S$ should give complex structures on the last tangent bundle. For example
if we find a natural section $s: T{\cal B}_S \hookrightarrow {\cal U}_{SBS}$ such that the image is perpendicular to the fibers with respect to the Kahler form, lifted from
$\mathbb{P} \Gamma(M, L)$ then the restriction $\tilde I|_{s(T {\cal B}_S)}$ must be an integrable complex structure.  We call this type of construction ``complexification''
of the moduli space ${\cal B}_S$: realizing this programme we get a complex space with ${\cal B}_S$ sitting inside as a completely real subspace.

  $$$$

{\bf References:}
        
[1] A. Gorodentsev, A. Tyurin, {\it ``Abelian lagrangian algebraic geometry''}, Izvestiya: Mathematics, 2001, 65:3, 437–467; 

[2] N. Tyurin, {\it ``Special Bohr - Sommerfeld lagrangian submanifolds''}, \\ Izvestiya: Mathematics, 2016, 80:6, 1257–1274;

[3] J. \'{S}niatycki, {\it ''Geometric Quantization and Quantum Mechanics``}, Appl. Math. Sci, 30, Springer, New York - Berlin, 1980.

\end{document}